\documentclass[11pt]{amsart}

\usepackage{amsmath,amssymb,amsthm,amsrefs,a4wide}
\usepackage{latexsym}
\usepackage{indentfirst}
\usepackage{graphicx}
\usepackage{placeins}
\usepackage{booktabs}
\usepackage{algorithm}
\usepackage{algorithmic}
\usepackage{multirow}
\usepackage{color}

\numberwithin{figure}{section}

\theoremstyle{plain}

\theoremstyle{definition}

\newtheorem{example}{Example}[section]

\theoremstyle{remark}

\DeclareMathOperator{\argmin}{arg min}

\newcommand{\diag}{\mathrm{diag}}

\newcommand{\bbm}{\begin{bmatrix}}
\newcommand{\ebm}{\end{bmatrix}}
\newcommand{\R}{\mathbb{R}}
\newcommand{\C}{\mathbb{C}}

\newcommand{\T}{\mathsf{T}}

\newcommand{\nx}{{n}}
\newcommand{\nz}{{n_z}}
\newcommand{\nc}{{n_c}}
\newcommand{\nl}{{n_l}}

\newcommand{\Bh}{\hat{B}}

\newcommand{\bb}{\mathbf{b}}
\newcommand{\bu}{\mathbf{u}}

\newcommand{\bbh}{\hat{\bb}}

\newcommand{\muA}{{\mu_A}}
\newcommand{\muB}{{\mu_B}}
\newcommand{\muC}{{\mu_C}}

\begin{document}

\title[Sparse free deconvolution under unknown noise level via eigenmatrix]{Sparse free deconvolution under unknown noise level via eigenmatrix}

\author[]{Lexing Ying} \address[Lexing Ying]{Department of Mathematics, Stanford University,
  Stanford, CA 94305} \email{lexing@stanford.edu}

\thanks{This work is partially supported by NSF grant DMS-2208163. The author thanks Xiucai Ding for helpful discussions.}

\keywords{Free deconvolution, sparse spectral measure, eigenmatrix.}

\subjclass[2010]{46L54, 62H12, 65R32.}

\begin{abstract}
  This note considers the spectral estimation problems of sparse spectral measures under unknown noise levels. The main technical tool is the eigenmatrix method for solving unstructured sparse recovery problems. When the noise level is determined, the free deconvolution reduces the problem to an unstructured sparse recovery problem to which the eigenmatrix method can be applied. To determine the unknown noise level, we propose an optimization problem based on the singular values of an intermediate matrix of the eigenmatrix method. Numerical results are provided for both the additive and multiplicative free deconvolutions.
\end{abstract}

\maketitle

\section{Introduction}\label{sec:intro}

This note considers free deconvolution problems of sparse spectral measures under unknown noise levels. In the additive setting, assume that $A$ is a real $N\times N$ symmetric matrix with an unknown sparse spectral measure $\mu_A$. Let $B$ be a Wigner matrix with an unknown noise level $\sigma$, i.e., the off-diagonal and diagonal entries are Gaussian with variance $\frac{\sigma^2}{N}$ and $\frac{2\sigma^2}{N}$, respectively. The spectral measure $\mu_B$ of $B$ follows the semicircle law with parameter $\sigma$ in the large $N$ limit. Given the spectral measure $\mu_C$ of $C=A+B$, the task is to recover $\sigma$ and $\mu_A$.

In the multiplicative setting, assume that $A$ is a real $N\times N$ symmetric positive definite matrix with an unknown sparse spectral measure $\mu_A$. Let $B$ be a Wishart matrix with an unknown dimension-to-sample size ratio $q$, i.e., $B$ is statistically equivalent to $\frac{1}{T} \sum_{t=1}^T X_t X_t^\T$ with $q=N/T$ with $X_t \sim \mathcal{N}(0, I_{N\times N})$. The spectral measure $\mu_B$ of $B$ follows the Marchenko-Pastur law with parameter $q$ in the large $N$  limit. Given the spectral measure $\mu_C$ of $C = \sqrt{A} B \sqrt{A}$, the task is to recover $q$ and $\mu_A$.

These two problems have many applications in statistics and data science. The first one, sometimes referred to as the deformed Wigner model, comes up when the data matrix $A$ is polluted by entrywise independent noise with unknown variance. The second problem appears in the estimation of the covariance matrix $A$, where $q$ is not known due to either the lack of information about $T$ or the dependence between the samples $\{X_t\}$.

\subsection{Related work.}
Spectral estimation has been an active field in the past two decades. For covariance matrix estimation (the multiplicative setting), many methods have been proposed over the years, including linear shrinkage \cite{ledoit2004well}, methods based on optimizing over the Marcenko-Pastur equation directly \cite{el2008spectrum}, the highly successful nonlinear shrinkage methods \cite{ledoit2011eigenvectors,ledoit2012nonlinear,ledoit2015spectrum,ledoit2020analytical}, and the moment-based method \cite{kong2017spectrum}. For the deformed Wigner model (the additive setting), there has also been a significant body of work, mostly for nonlinear shrinkage methods, including \cite{gavish2014optimal,bun2016rotational,lolas2021shrinkage}.

There is little work on using the free deconvolution directly. In \cite{arizmendi2020subordination}, the subordination method reduces the task to a classical deconvolution problem, which is then solved via Tikhonov regularization with carefully chosen regularization parameters.

In most of the work mentioned above, the noise level $\sigma$ or $q$ is given. The case with unknown noise level is much less explored.

\subsection{Contributions.} In this note, we consider the more challenging case where the noise level $\sigma$ or $q$ is unknown. Without prior information on the spectral measure $\mu_A$ of $A$, the problems are ill-defined. In order to work with a well-defined setting, we assume that $\mu_A$ is sparsely supported.

The main technical tool is the recently proposed eigenmatrix method \cite{ying2024eigenmatrix,ying2025multidimensional} for solving the unstructured sparse recovery problem. The word "unstructured" refers to the fact that the sample locations can be arbitrary, as will be the case for the deconvolution problems considered here.

Assume for a moment that the noise level is already determined. By using the R-transform and S-transform \cite{voiculescu1986addition,voiculescu1987multiplication} from free probability \cite{mingo2017free,potters2020first}, the spectral estimation problem can be cast into a classical inverse problem of the Cauchy integral with observations at unstructured locations. Since $\mu_A$ is sparse, this can be solved directly with the eigenmatrix method.

To determine the noise level, we propose a novel optimization problem in which the loss function is based on the singular value of an intermediate matrix used in the eigenmatrix method. Optimizing this objective function gives an accurate estimate of the noise level.

The rest of the note is organized as follows. Section \ref{sec:em} reviews the eigenmatrix method. Section \ref{sec:add} describes the additive case. Section \ref{sec:mul} considers the multiplicative case. Section \ref{sec:disc} concludes with a discussion for future work.

\section{Eigenmatrix} \label{sec:em}

This section provides a short review of the eigenmatrix method for unstructured sparse recovery problems. Let $X$ be the parameter space. $G(z,x)$ is a kernel function defined for $x\in X$ at sample $z$, and is assumed to be analytic in $x$. In this note, $X$ is an interval of $\R$, $z\in \C$, and $G(z,x)=\frac{1}{z-x}$. Suppose that
\begin{equation}
  f(x) = \sum_{k=1}^{\nx} w_k \delta(x-x_k)
  \label{eq:f}
\end{equation}
is an unknown sparse signal, where $\{x_k\}_{1\le k \le \nx}$ are the spike locations and $\{w_k\}_{1\le k \le \nx}$ are the spike weights. The observed quantity is defined via the function
\begin{equation}
  u(z) :=\int_X G(z,x) f(x) dx = \sum_{k=1}^{\nx} G(z,x_k) w_k.
  \label{eq:u}
\end{equation}
Let $\{z_j\}_{1\le j \le n_z}$ be a set of $\nz$ unstructured samples. Suppose that we are given the noisy observations $u_j \approx u(z_j)$. The task is to recover the spikes $\{x_k\}$ and weights $\{w_k\}$ from $\{u_j\}$.

Define for each $x$ the column vector
\[
\bb_x:= 
\bbm
G(z_j,x)
\ebm_{1\le j \le \nz}
\]
in $\C^{n_z}$. Notice that $\bb_x$ is analytic in terms of $x$.

The first step is to construct a matrix $M \in \C^{n_z\times n_z}$ such that $M \bb_{x} \approx x \bb_{x}$ for $x\in X$, i.e., $M$ is the matrix with $(x,\bb_{x})$ as approximate eigenpairs for $x\in X$.  Numerically, it is more robust to use the normalized vector $\bbh_{x} = \bb_{x}/\|\bb_{x}\|$ since the norm of $\bb_{x}$ can vary significantly depending on $x$. The condition then becomes
\[
M \bbh_{x} \approx x \bbh_{x}, \quad x\in X.
\]
To construct $M$ numerically, we choose a Chebyshev grid $\{c_t\}_{1\le t \le \nc}$ of size $\nc$ on the interval $X$, where $\nc$ is sufficiently large yet the vectors $\{\bbh_{c_t}\}$ are numerically linearly independent. We then enforce this condition on this grid, i.e.,
\[
M \bbh_{c_t} \approx c_t \bbh_{c_t}.
\]
By defining the $\nz\times \nc$ matrix
\[
\Bh = 
\bbm
\bbh_{c_1} & \ldots & \bbh_{c_\nc}
\ebm
\]
and the $\nc\times \nc$ diagonal matrix $\Lambda = \diag(c_t)$, the previous condition can be written in a matrix form as
\[
M \Bh \approx \Bh \Lambda.
\]
Because the columns of $\Bh$ are numerically linearly independent, we set the eigenmatrix as
\begin{equation}
  M := \Bh \Lambda \Bh^+, \label{eq:M}
\end{equation}
where the pseudoinverse $\Bh^+$ is computed by thresholding the singular values of $\Bh$. In practice, the thresholding value is chosen so that the norm of $M$ is bounded by a small constant. 

Once $M$ is ready, the rest follows, for example, the ESPRIT algorithm \cite{roy1989esprit}. Define the vector
\[
\bu=
\bbm   
u_j
\ebm_{1\le j \le \nz}
\]
in $\C^{n_z}$, where $u_j$ are the noisy observations. Notice that $\bu \approx \sum_k \bb_{x_k} w_k$. Consider the matrix
\begin{equation}
  T \equiv \bbm
  \bu & M\bu & \ldots & M^\nl \bu
  \ebm
  \label{eq:T}
\end{equation}
with $\nl > \nx$, obtained from applying $M$ repeatitively to $\bu$. Since $\bu \approx \sum_k \bb_{x_k} w_k$ and $M\bb_{x}\approx x \bb_{x}$,
\[
T = \bbm  \bu & M\bu & \ldots & M^\nl \bu \ebm
\approx
\bbm \bb_{x_1} & \ldots & \bb_{x_\nx} \ebm
\bbm w_1 & &\\& \ddots & \\& & w_\nx\ebm
\bbm
1 &  x_1 & \ldots & (x_1)^\nl \\
\vdots & \vdots & \ddots & \vdots \\
1 &  x_\nx & \ldots & (x_\nx)^\nl
\ebm.
\]
Let $U S V^*$ be the rank-$\nx$ truncated SVD of $T$. The matrix ${V}^*$ then satisfies
\[
{V}^* \approx P
\bbm
1 &  x_1 & \ldots & (x_1)^\nl \\
\vdots & \vdots & \ddots & \vdots \\
1 &  x_\nx & \ldots & (x_\nx)^\nl
\ebm,
\]
where $P$ is an unknown non-degenerate $\nx \times \nx$ matrix. Let $Z_L$ and $Z_H$ be the submatrices obtained by excluding the last column and the first column of $V^*$, respectively, i.e.,
\[
Z_L \approx P
\bbm
1      & \ldots & (x_1)^{\nl-1} \\
\vdots & \ddots & \vdots \\
1      & \ldots & (x_\nx)^{\nl-1}
\ebm,
\quad
Z_H \approx P
\bbm
x_1 & \ldots & (x_1)^\nl \\
\vdots & \ddots & \vdots \\
x_\nx & \ldots & (x_\nx)^\nl
\ebm.
\]
By forming $Z_H (Z_L)^+$ and noticing
\[
Z_H (Z_L)^+ \approx
P
\bbm
x_1 &  & \\
& \ddots & \\
& & x_\nx
\ebm
P^{-1},
\]
one obtains the estimates for $\{x_k\}$ by computing the eigenvalues of $Z_H (Z_L)^+$.

With the estimates for $\{x_k\}$ available, the least square solution of 
\[
\min_{w_k} \sum_j \left|\sum_k G(s_j,x_k) w_k - u_j\right|^2 
\]
gives the estimators for $\{w_k\}$.

\section{Additive deconvolution}\label{sec:add}

To address the additive case, we leverage the R-transform \cite{voiculescu1986addition}. Given a spectral measure $\mu$, the Cauchy integral defines a correspondence between $z$ and $g = \int \frac{1}{z-x} d\mu(x)$. The map from $z$ to $g$ is the Stieltjes transform, denoted by $g_\mu(z)$. Its inverse map from $g$ to $z$ is well-defined for sufficiently small values of $g$ and is denoted by $z_\mu(g)$. The R-transform is then defined as 
\[
r_\mu(g) = z_\mu(g)-\frac{1}{g}.
\]

Since $C=A+B$, in the large dimension limit $\mu_C = \mu_A \boxplus \mu_B$ and  
\[
r_\muC(g) = r_\muA(g) + r_\muB(g).
\]

\subsection{Known $\sigma$}\label{sec:knownsigma}

Assume for now that $\sigma$ is determined. Then $r_\muB(g) = \sigma^2 g$ from the semicircle law. Due to its sparsity, $\mu_A = \sum_k \delta_{x_k} w_k$. The task is to recover $x_k$ and $w_k$.

Given $\muC$, we choose $\{z_j\}$ to be a set of points on an ellipsis around $\mu_C$ and compute
\[
g_j = \int \frac{1}{z_j - x} d \mu_C(x), 
\]
i.e., $g_\muC(z_j) = g_j$ and $z_\muC(g_j) = z_j$.

From $r_\muC(g) = r_\muA(g) + r_\muB(g)$ and $r_\muB(g) = \sigma^2 g$, we have
\[
r_\muA(g_j) = r_\muC(g_j) - \sigma^2 g_j  \quad     z_\muA(g_j) = z_\muC(g_j) - \sigma^2 g_j = z_j - \sigma^2 g_j.
\]
Define $z_j' = z_j - \sigma^2 g_j$. Then $z_j' = z_\muA(g_j)$ and $g_\muA(z_j') = g_j$, i.e., $(z_j', g_j)$ are samples for the Stieltjes transform of $\mu_A$:
\[
g_j = \int \frac{1}{z_j'-x} d \mu_A (x) = \sum_k \frac{1}{z_j'-x_k} w_k.
\]
Since the locations $\{z_j'\}$ are not a priori controlled, recovering $\{x_k\}$ and $\{w_k\}$ is a sparse, unstructured recovery problem.

Next, we apply the eigenmatrix method. More specifically, set $X$ to be the shortest interval that covers $\muC$. Treat $\{z_j'\}$ as the samples and $\{g_j\}$ as the observed data.  Define $\bb_x = \bbm \frac{1}{z_j'-x} \ebm_j$. Choose a Chebyshev grid $\{c_t\}$ and construct $M$ such that $M \bb_{c_t} \approx c_t \bb_{c_t}$. Then, define $\bu = \bbm g_j\ebm_j$ and form the $T$ matrix \eqref{eq:T} to recover $\{x_k\}$ and $\{w_k\}$.

\subsection{Unknown $\sigma$}

Let us consider how to determine $\sigma$. For a value of $\sigma$, since the $T$ matrix is a function of $\sigma$, we denote it by $T(\sigma)$. The main observation is the following: at the correct $\sigma$ value, $T(\sigma)$ has a numerical rank equal to $\nx$, while for other values of $\sigma$, the rank is higher.

\begin{figure}[h!]
  \centering
  \includegraphics[scale=0.28]{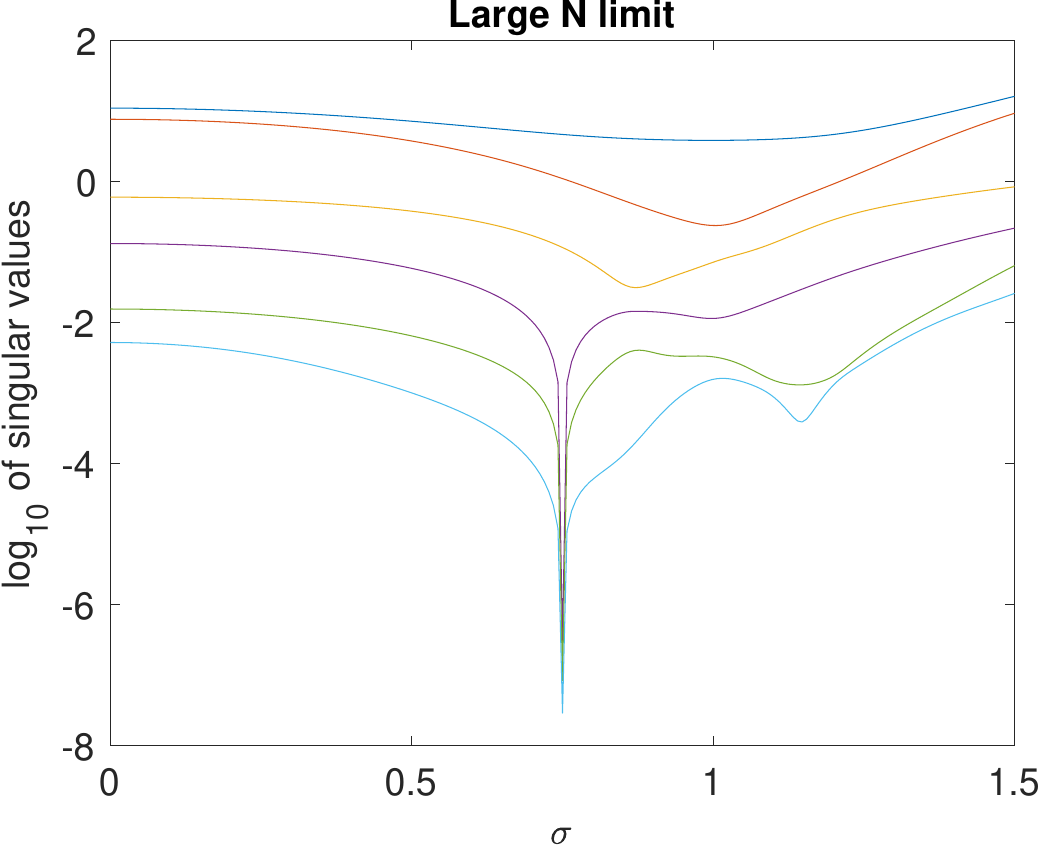}
  \includegraphics[scale=0.28]{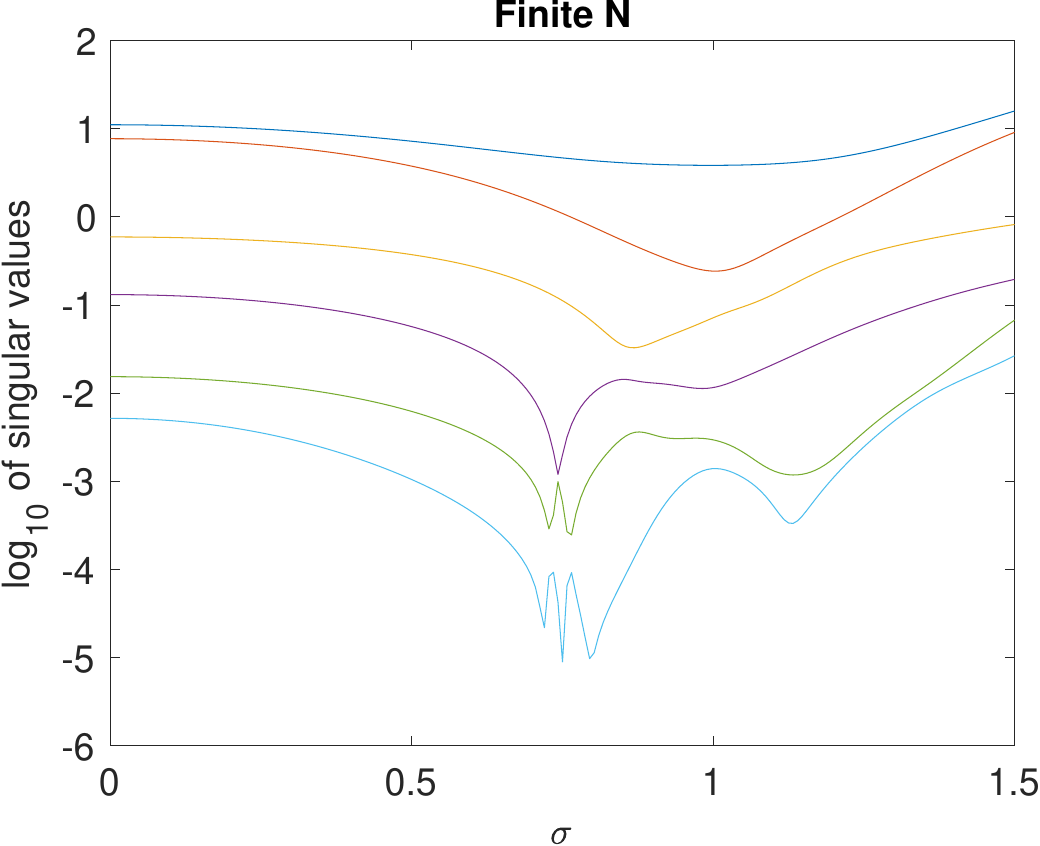}
  \includegraphics[scale=0.28]{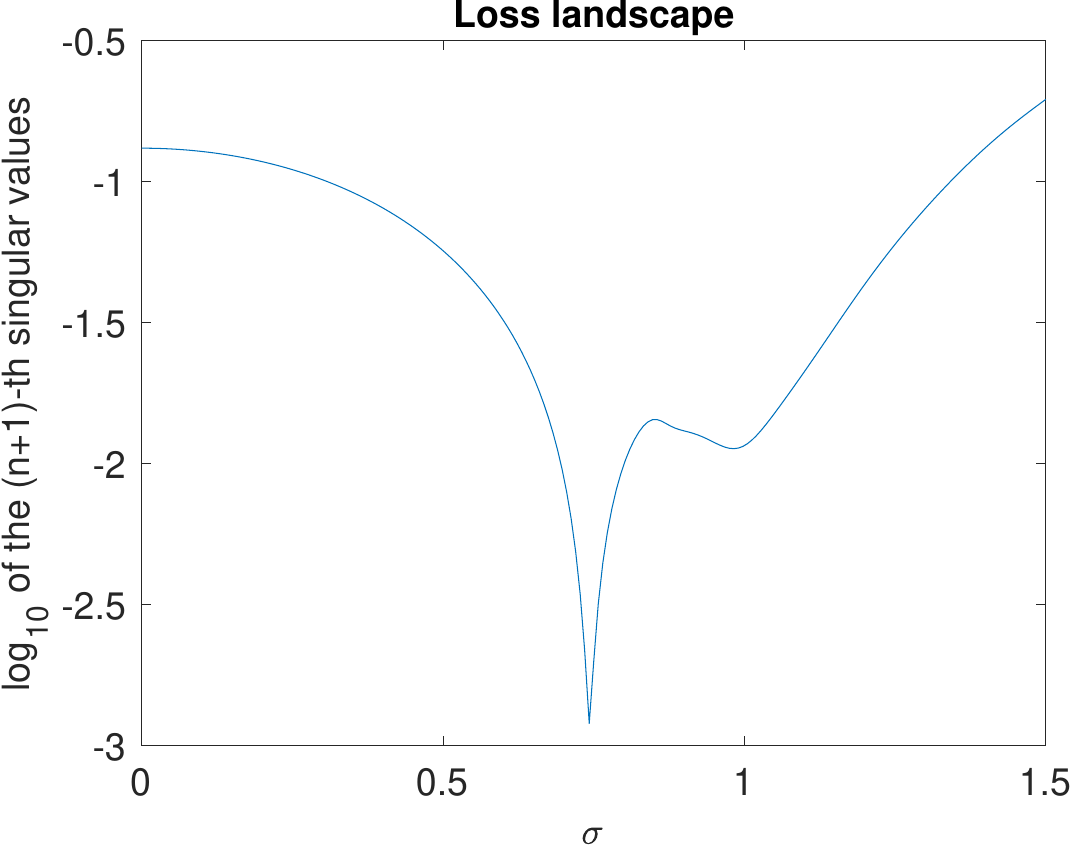}
  \caption{Logarithm of singular values of $T(\sigma)$ as a function of $\sigma$. Left: large $N$ limit. Middle: $N=1024$. Right: loss landscape.}
  \label{fig:aloss}
\end{figure}

Consider an example where $\sigma=0.75$ and $\muA$ has $n=3$ spikes.  Figure \ref{fig:aloss} plots the singular values of $T(\sigma)$ in logarithmic scale as a function of $\sigma$, i.e., the $i$-th curve from the top stands for the $i$-th singular value. The left plot shows the large $N$ limit. At $\sigma=0.75$, starting from the 4th curve, the singular value drops numerically to zero. This plot indicates that the right $\sigma$ and $n$ can be easily detected in the large $N$ limit.

The middle plot shows the situation at $N=1024$. The overall trends remain the same. However, due to the finiteness of $N$, it is harder to find an automatic procedure to identify $\sigma$ and $n$ together. In what follows, we assume $\nx$ is given. Here, the ($n$+1)-th (i.e., 4th) singular value curve clearly shows a global minimum at the right $\sigma$.

To automate this process, we formulate it as a minimization problem
\[
\hat{\sigma} = \argmin_\sigma \log(s_{\nx+1}(T(\sigma))),
\]
where $s_{\nx+1}(\cdot)$ refers to the ($n$+1)-th singular value. The right plot shows that this loss landscape can be quite non-convex, which can cause local optimization to get stuck. Since this is a one-dimensional optimization problem, we adopt a two-step procedure. The first step performs a grid search on a coarse grid and gives a good initial guess $\sigma_{\text{init}}$. Starting from $\sigma_{\text{init}}$, the second step performs a local search to identify the optimal $\hat{\sigma}$. Finally, given $\hat{\sigma}$, we can compute the $(z_j', g_j)$ pairs of $\mu_A$ as in Section \ref{sec:knownsigma} to recover $\{x_k\}$ and $\{w_k\}$ via eigenmatrix.

Below, we present a few numerical examples.

\begin{example} The parameters are
  \begin{itemize}
  \item $N=1024$.
  \item $\sigma = 0.25$.
  \item $\mu_A = \frac{1}{4} \delta_{-1} + \frac{1}{2} \delta_{0.2} + \frac{1}{4} \delta_{1}$.
  \end{itemize}
  
  \begin{figure}[h!]
    \centering
    \includegraphics[scale=0.3]{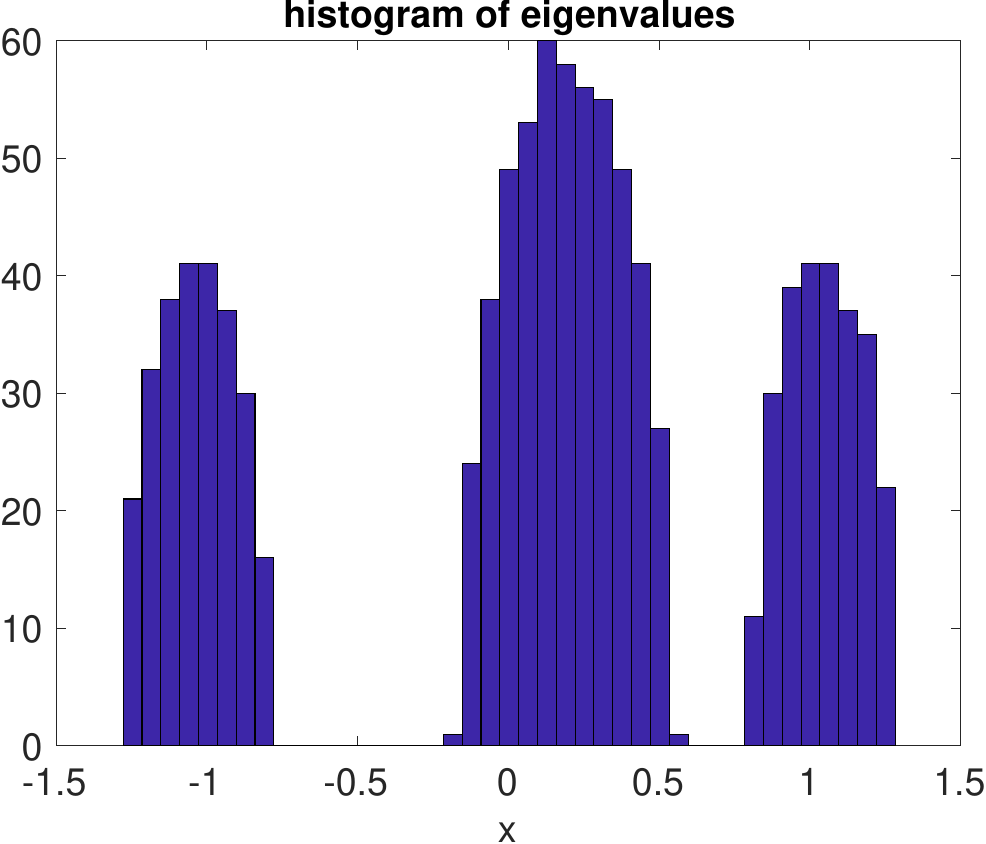}
    \includegraphics[scale=0.3]{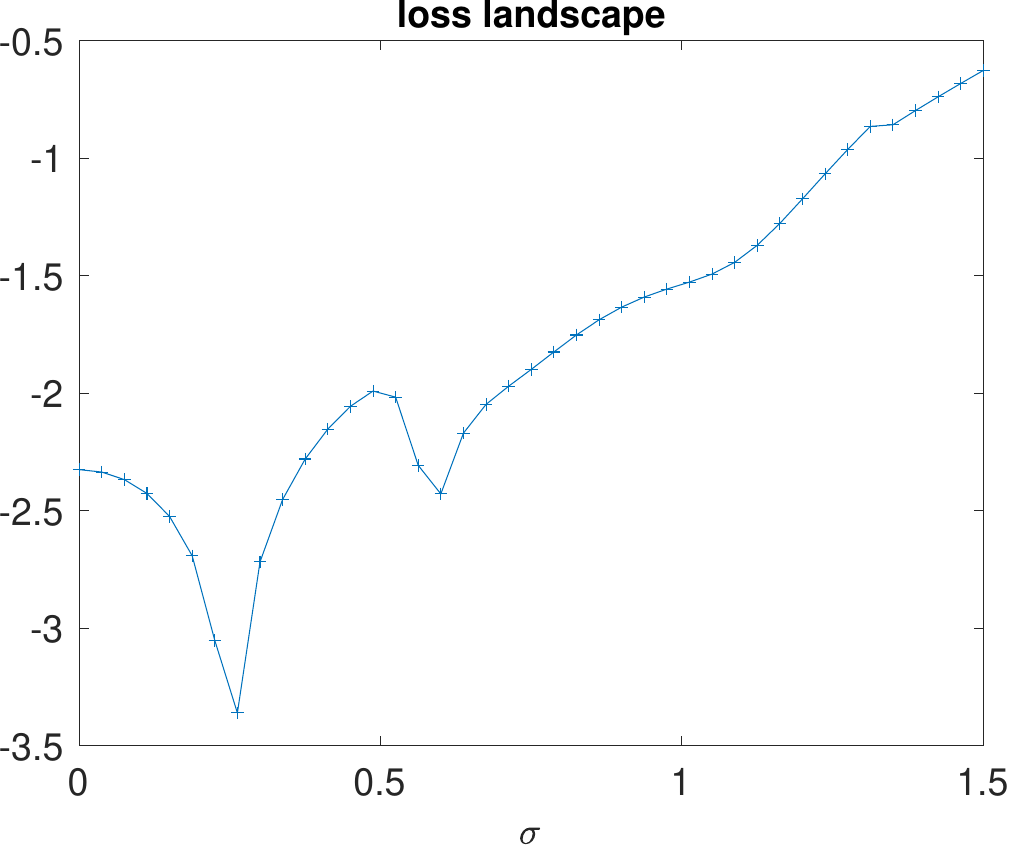}
    \includegraphics[scale=0.3]{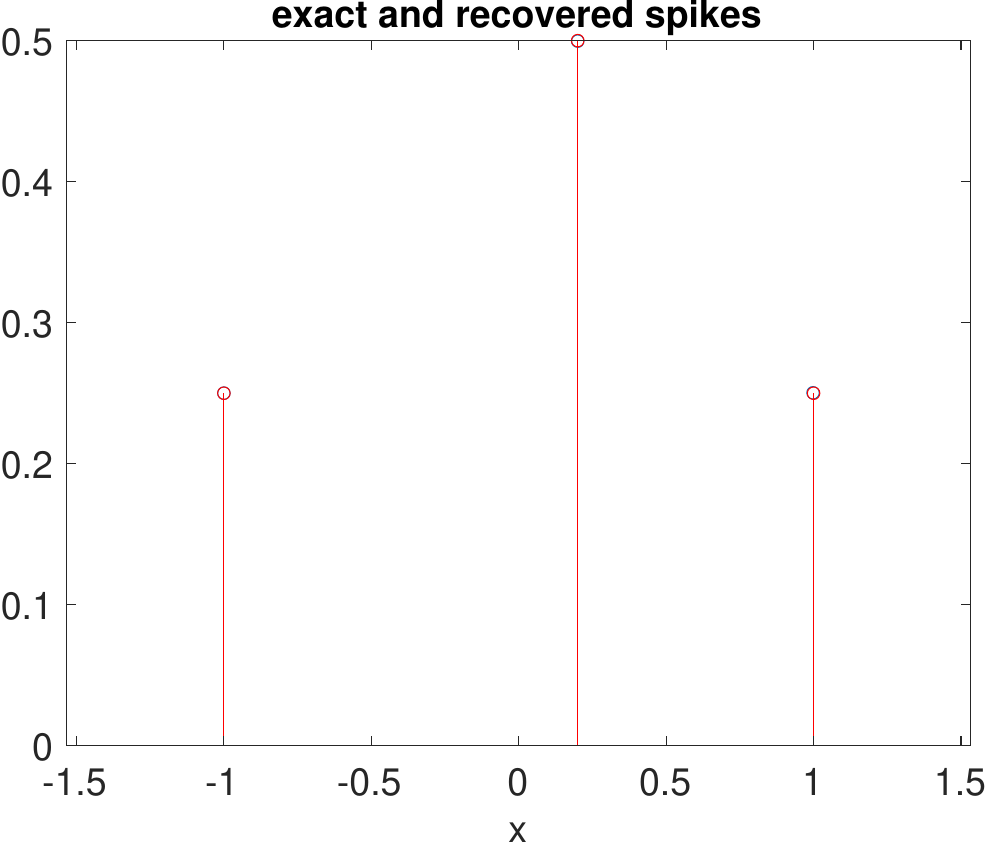}
    \caption{Left: the histogram of the eigenvalues of $C$. Middle: the loss landscape as a function of $\sigma$, and the global minimum is around the ground truth $\sigma$.  Right: the exact (red) and reconstructed (blue) $\muA$.}
    \label{fig:a1}
  \end{figure}
  
  Figure \ref{fig:a1} summarizes the results. The first plot gives the histogram of the eigenvalues of $C$. The second plot shows the loss landscape as a function of $\sigma$, and the global minimum is around the ground truth $\sigma$. The last plot shows the exact (red) and reconstructed (blue) spectral measure of $A$.
\end{example}

\begin{example} The parameters are
  \begin{itemize}
  \item $N=1024$.
  \item $\sigma = 0.75$.
  \item $\mu_A = \frac{1}{4} \delta_{-1} + \frac{1}{2} \delta_{0.2} + \frac{1}{4} \delta_{1}$.
  \end{itemize}
  
  \begin{figure}[h!]
    \centering
    \includegraphics[scale=0.3]{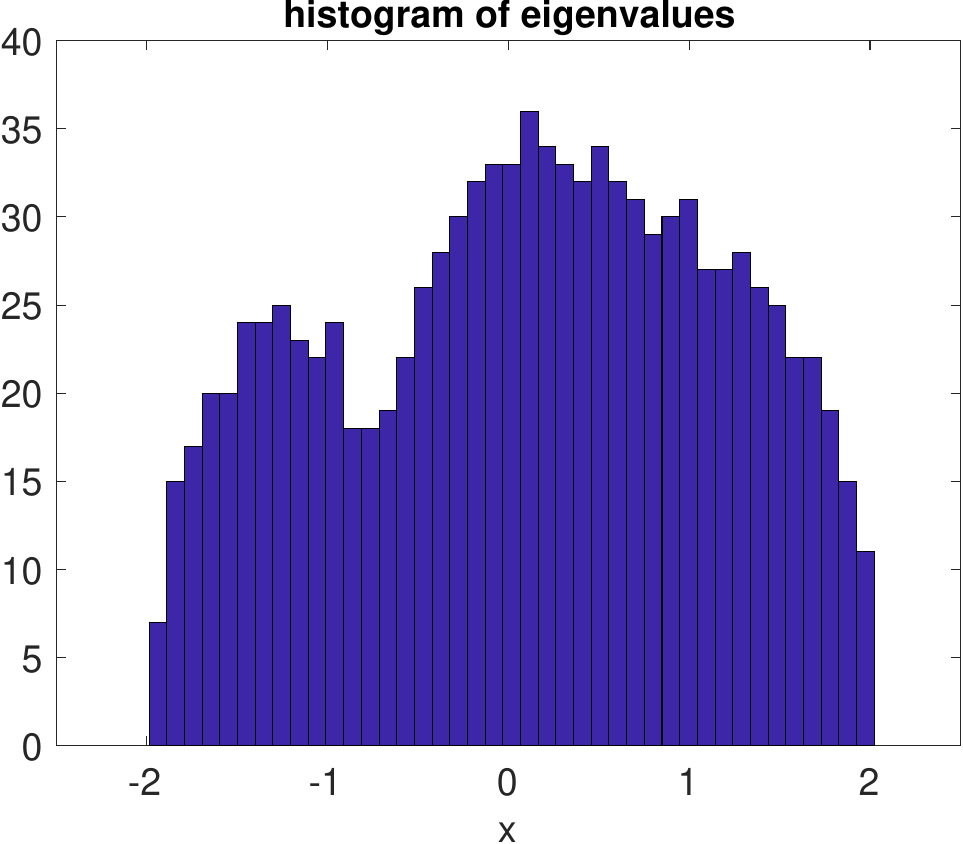}
    \includegraphics[scale=0.3]{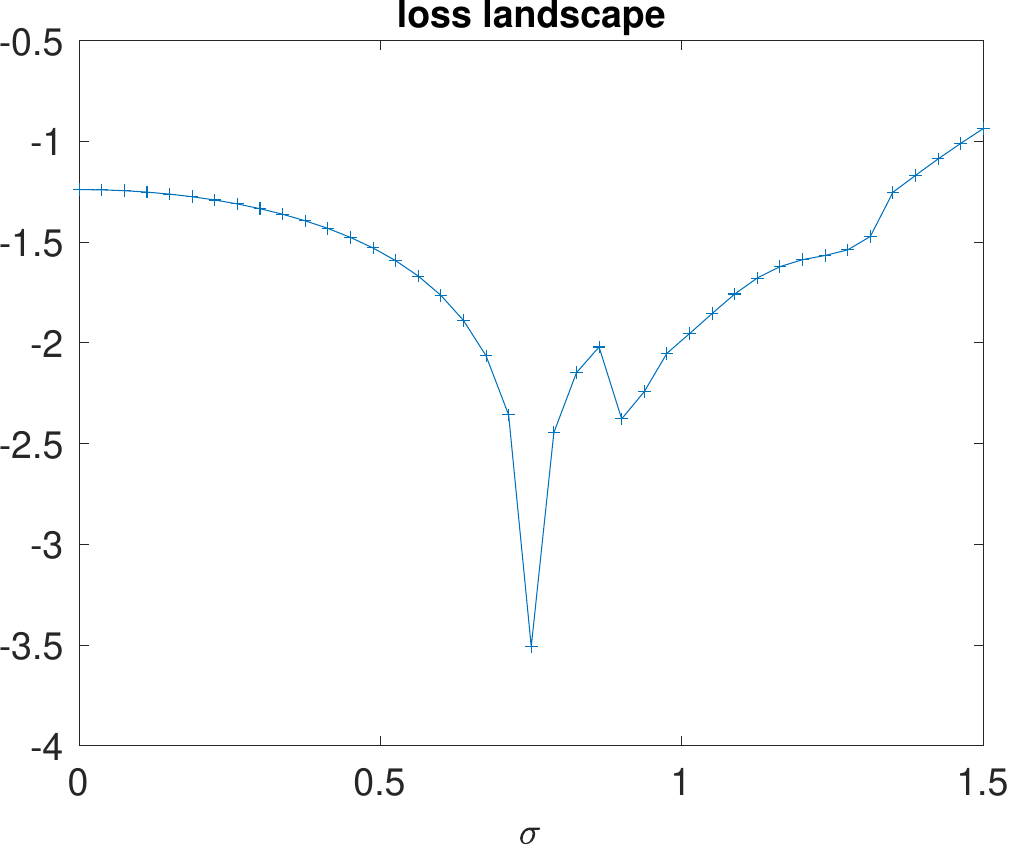}
    \includegraphics[scale=0.3]{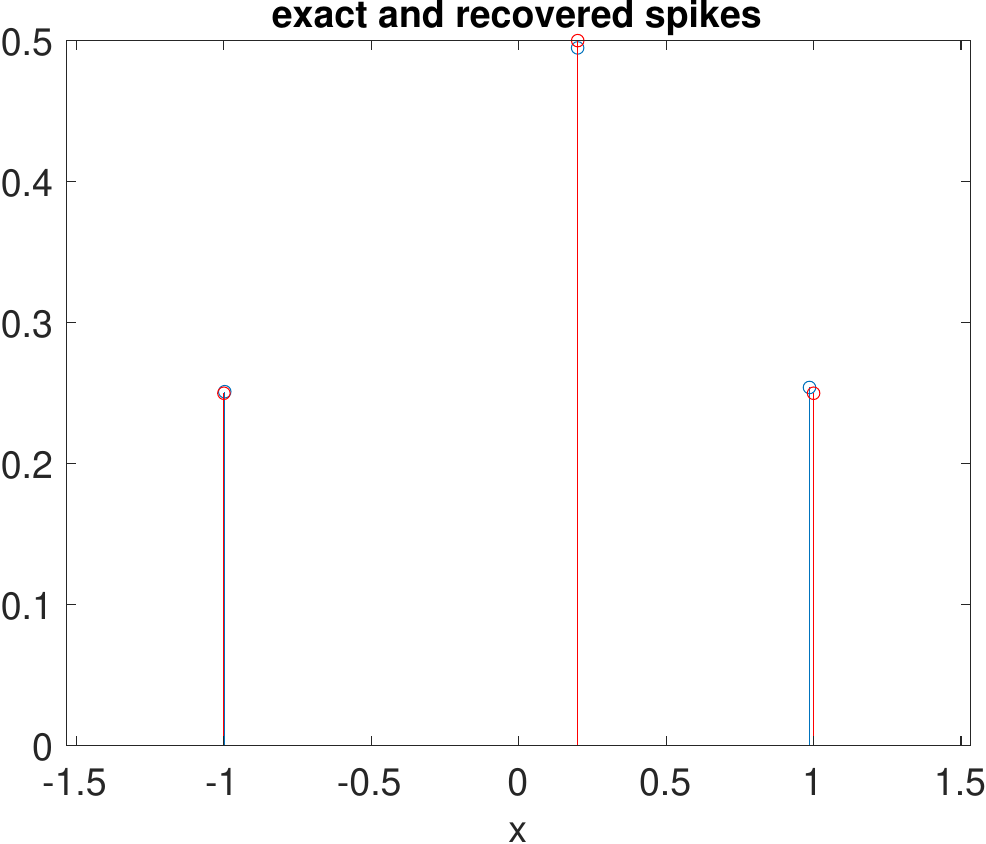}
    \caption{Left: the histogram of the eigenvalues of $C$. Middle: the loss landscape as a function of $\sigma$, and the global minimum is around the ground truth $\sigma$.  Right: the exact (red) and reconstructed (blue) $\muA$.}
    \label{fig:a2}
  \end{figure}
  
  Figure \ref{fig:a2} summarizes the results. The meanings of the plots are the same as in the previous example.


\end{example}

\begin{example} The parameters are
  \begin{itemize}
  \item  $N=1024$.
  \item  $\sigma = 1.25$.
  \item  $\mu_A = \frac{1}{4} \delta_{-1} + \frac{1}{2} \delta_{0.2} + \frac{1}{4} \delta_{1}$.
  \end{itemize}
  
  \begin{figure}[h!]
    \centering
    \includegraphics[scale=0.3]{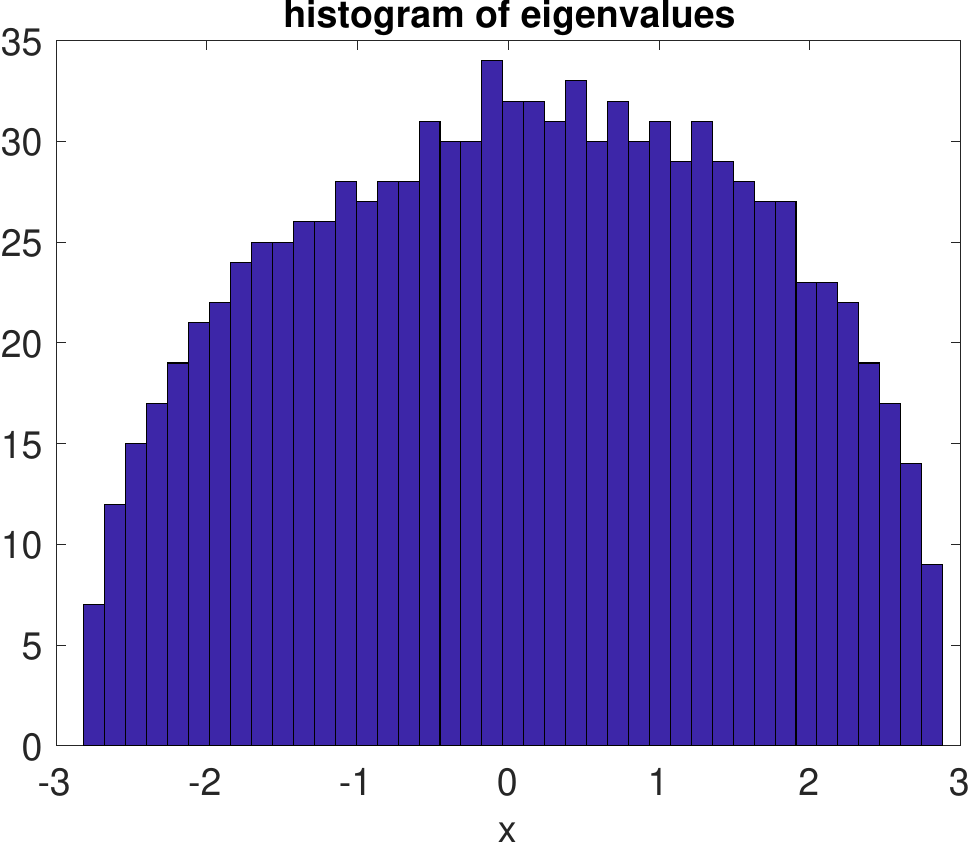}
    \includegraphics[scale=0.3]{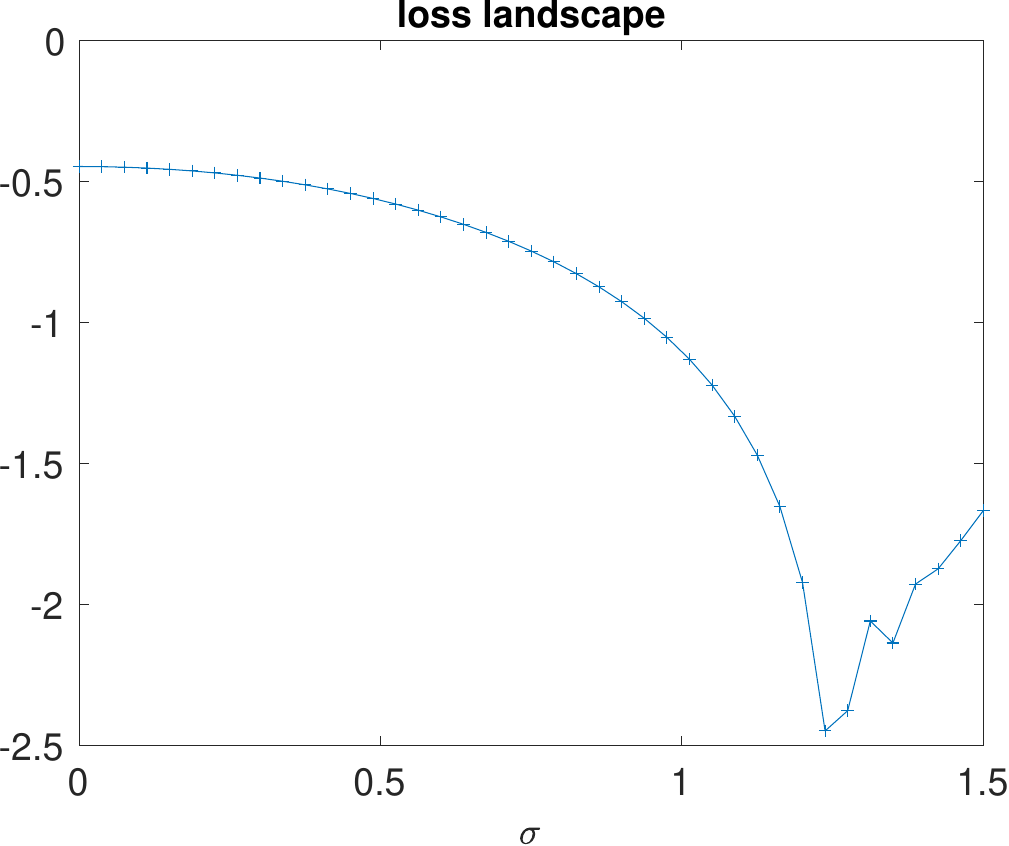}
    \includegraphics[scale=0.3]{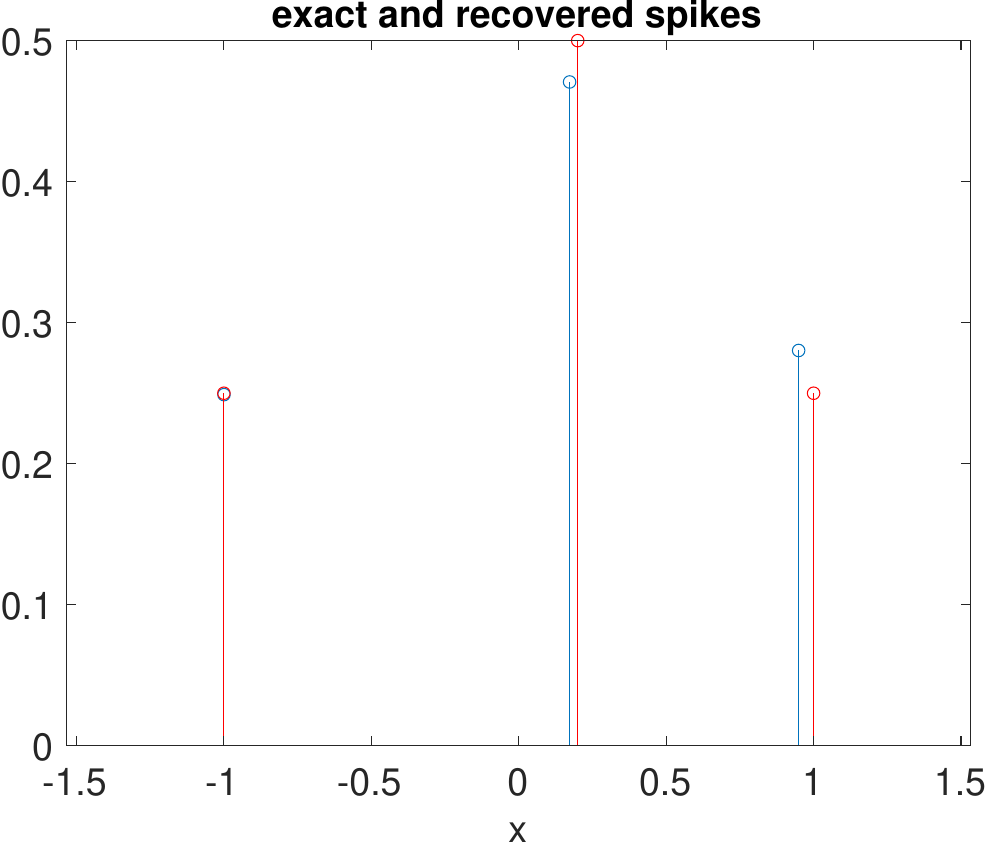}
    \caption{Left: the histogram of the eigenvalues of $C$. Middle: the loss landscape as a function of $\sigma$, and the global minimum is around the ground truth $\sigma$.  Right: the exact (red) and reconstructed (blue) $\muA$.}
    \label{fig:a3}
  \end{figure}

  Figure \ref{fig:a3} summarizes the results. The meanings of the plots are the same as in the previous example.

  
\end{example}

\section{Multiplicative deconvolution}\label{sec:mul}

To address the multiplicative case, we leverage the S-transform \cite{voiculescu1987multiplication}. Given a spectral measure $\mu$, the Cauchy integral establishes a map between $z$ and $g = \int \frac{1}{z-x} d\mu(x)$. In addition, introduce the corresponding $t = z g-1$ and $s = \frac{t+1}{t z}$. The map from $t$ to $s$, well-defined for sufficiently small values of $t$, is the S-transform, denoted by $s_\mu(t)$.

Since $C = \sqrt{A} B \sqrt{A}$, in the large dimension limit $\mu_C = \mu_A \boxtimes \mu_B$ and 
\[
s_\muC(t) = s_\muA(t) s_\muB(t).
\]

\subsection{Known $q$}\label{sec:knownq}

Assume for now that $q$ is determined. $s_\muB(t) = \frac{1}{1+qt}$ from the Marchenko-Pastur law. Due to its sparsity, $\mu_A = \sum_k \delta_{x_k} w_k$. The task is to recover $\{x_k\}$ and $\{w_k\}$.

Given $\muC$, we choose $\{z_j\}$ to be a set of points on an ellipsis around $\mu_C$, compute
\[
g_j = \int \frac{1}{z_j - x} d \mu_C(x)
\]
and further define
\[
t_j = z_j g_j -1, \quad s_j  = \frac{t_j+1}{t_j z_j},
\]
i.e. $s_\muC(t_j) = s_j$.

From $s_\muC(t) = s_\muA(t) s_\muB(t)$ and $s_\muB(t) = \frac{1}{1+qt}$, we have 
\[
s_A(t_j) = s_C(t_j)/s_B(t_j) = s_C(t_j) (1+q t_j) = s_j (1+q t_j).
\]
Introduce $s_j' = s_j (1+q t_j)$. Then $(t_j, s_j')$ are samples of $s_\muA(t)$. i.e., $s_A(t_j)=s_j'$. Now further define
\[
z_j' = \frac{t_j+1}{t_j s_j'}, \quad g_j' = \frac{t_j+1}{z_j'}.
\]
The resulting $(z_j',g_j')$ are samples of the Stieltjes transform of $\mu_A$, i.e., 
\[
g_j' = \int \frac{1}{z_j'-x} d \mu_A (x) = \sum_k \frac{1}{z_j'-x_k} w_k.
\]
Since the locations $\{z_j'\}$ are not a priori controlled, recovering $\{x_k\}$ and $\{w_k\}$ is a sparse, unstructured recovery problem.

Next, we apply the eigenmatrix method. First, set $X$ to be the shortest interval that covers the spectrum of $C$. Treat $\{z_j'\}$ as the samples and $\{g_j'\}$ as the observed data.  Define $\bb_x = \bbm \frac{1}{z_j'-x} \ebm_j$.  Choose a Chebyshev grid $\{c_t\}$ and construct $M$ such that $M \bb_{c_t} \approx c_t \bb_{c_t}$.
Finally, define $\bu = \bbm g_j'\ebm_j$ and form the $T$ matrix to recover $\{x_k\}$ and $\{w_k\}$.

\subsection{Unknown $q$}

Let us discuss now how to determine $q$.  For a value of $q$, the $T$ matrix, as a function of $q$, will be denoted by $T(q)$. The main observation is similar to the additive case: for the correct $q$ value, $T(q)$ has a numerical rank equal to $\nx$, while for other values of $q$, the rank is higher.

\begin{figure}[h!]
  \centering
  \includegraphics[scale=0.28]{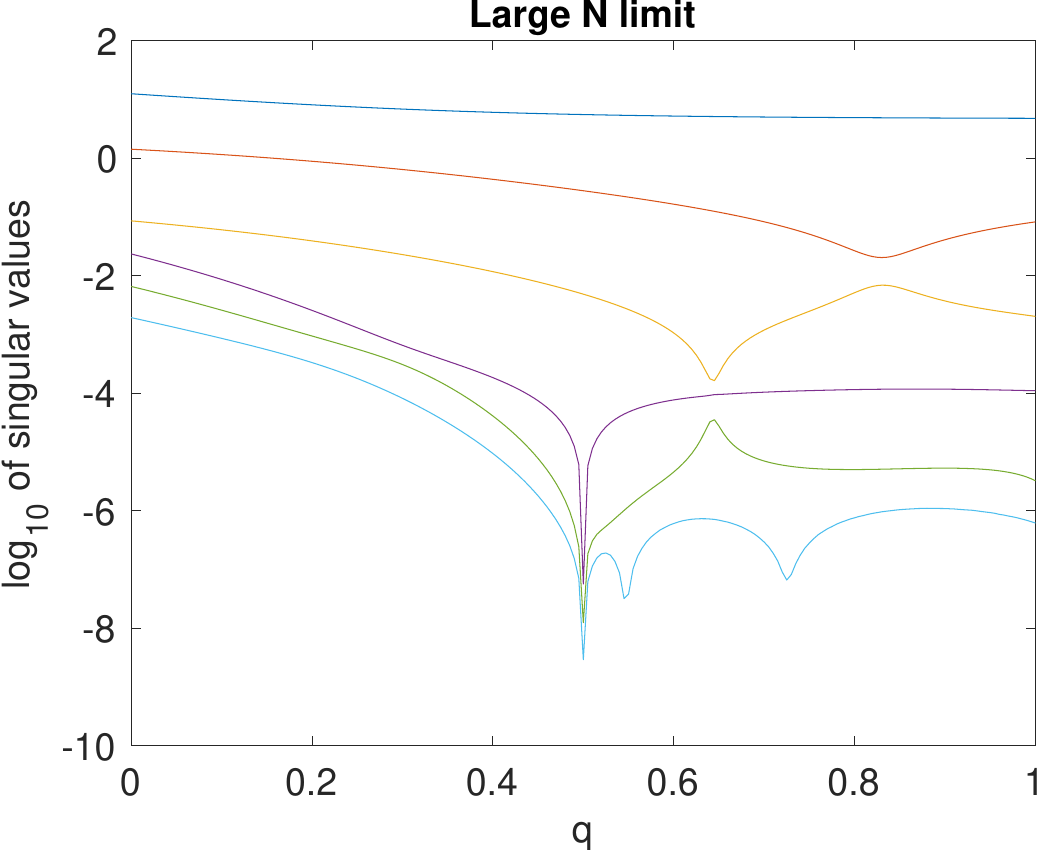}
  \includegraphics[scale=0.28]{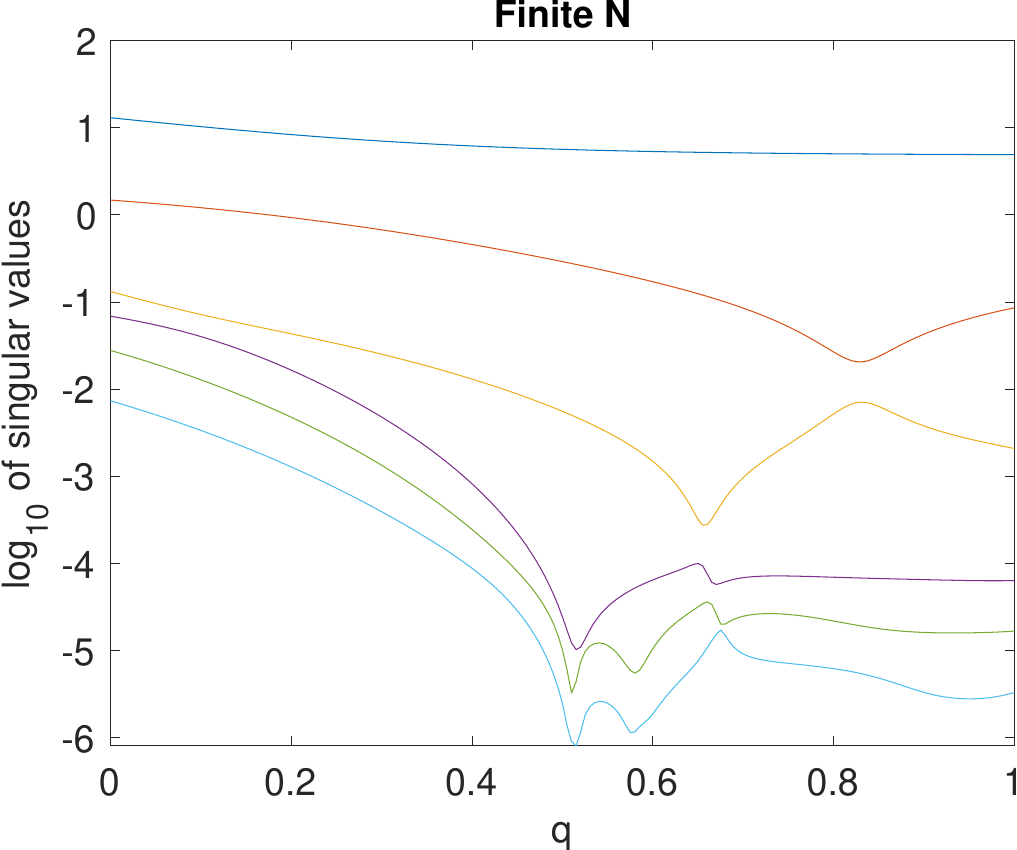}
  \includegraphics[scale=0.28]{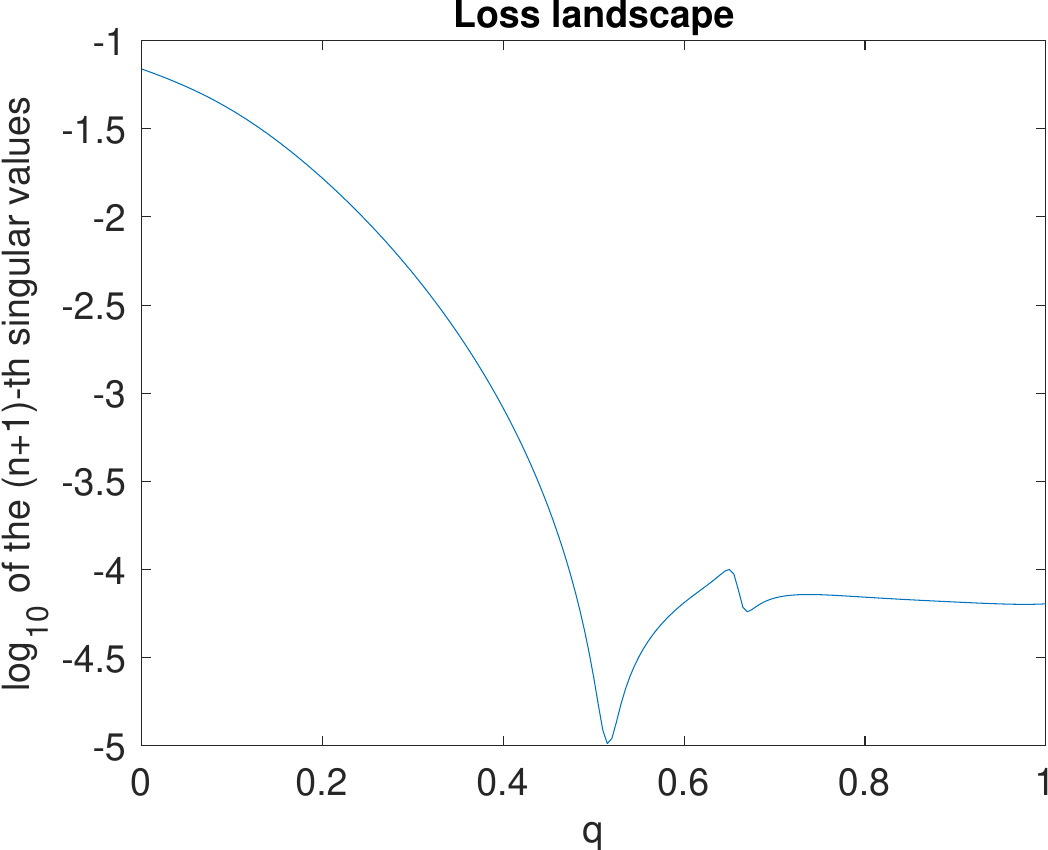}
  \caption{Logarithm of singular values of $T(q)$ as a function of $q$. Left: large $N$ limit. Middle: $N=1024$. Right: loss landscape.}
  \label{fig:mloss}
\end{figure}

Consider an example where $\muA$ has $n=3$ spikes and $q=0.5$. Figure \ref{fig:mloss} plots the singular values of $T(q)$ in logarithmic scale as a function of $q$. The left plot shows the large $N$ limit. At $q=0.5$, starting from the 4th curve, the singular value drops numerically to zero. This plot indicates the right $q$ and $n$ can be easily detected in the large $N$ limit.

The middle plot shows the case $N=1024$. The overall trends remain the same. However, due to the finiteness of $N$, it is harder to find an automatic procedure to identify $q$ and $n$ at the same time. In what follows, we assume $\nx$ is given. In this case, the ($n$+1)-th (i.e., 4th)  singular value curve clearly has a global minimum at the right $q$.

To automate this process, we formulate it again as a minimization problem
\[
\hat{q} = \argmin_q \log(s_{\nx+1}(T(q))),
\]
where $s_{\nx+1}(\cdot)$ refers to the ($n$+1)-th singular value. The right plot highlights that this loss landscape is not convex, so local optimization can get stuck. Since this is still a one-dimensional optimization problem, we again take a two-step procedure similar to the additive case.  First, we perform a coarse grid search to obtain a good initial guess $q_{\text{init}}$. Second, starting from $q_{\text{init}}$, a local search refines to the optimal $\hat{q}$. Finally, given $\hat{q}$, we can compute the $(z_j', g_j')$ pairs of $\mu_A$ as in Section \ref{sec:knownq} to find $\{x_k\}$ and $\{w_k\}$.

Below, we give a few numerical examples.

\begin{example} The parameters are
  \begin{itemize}
  \item $N=1024$.
  \item $q = 0.25$.
  \item $\mu_A = \frac{1}{3} \delta_{0.2} + \frac{1}{3} \delta_{0.6} + \frac{1}{3} \delta_{1}$.
  \end{itemize}
  
  \begin{figure}[h!]
    \centering
    \includegraphics[scale=0.27]{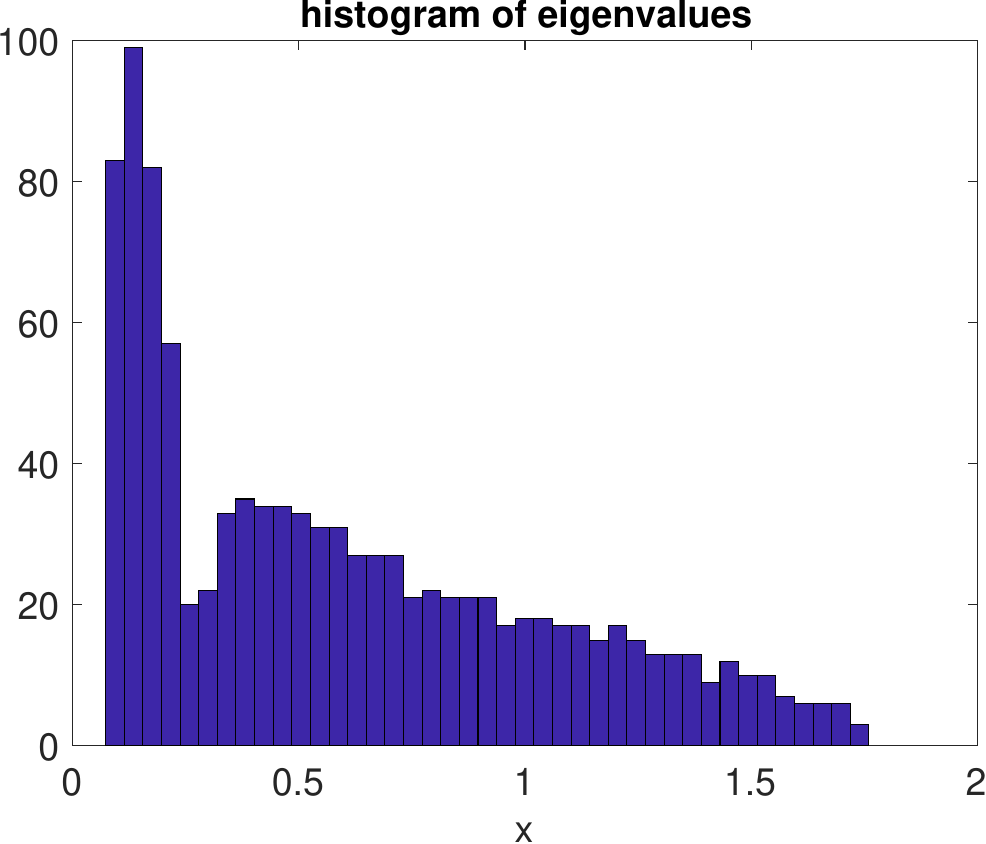}
    \includegraphics[scale=0.27]{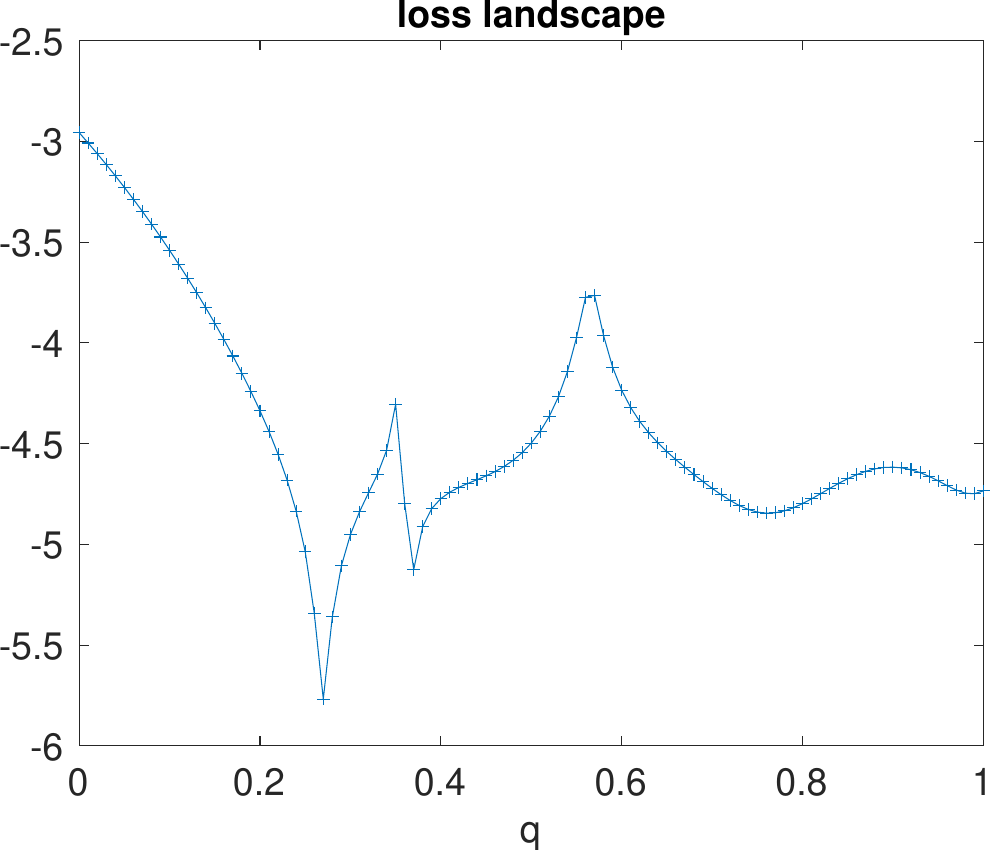}
    \includegraphics[scale=0.27]{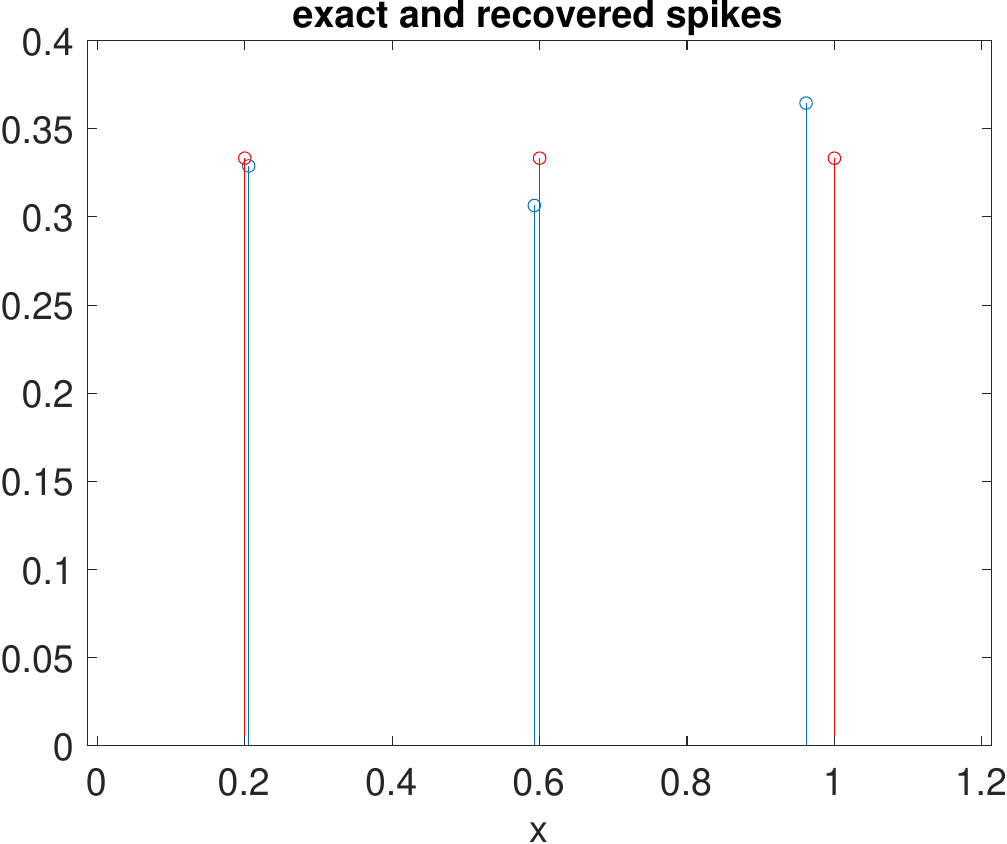}
    \caption{Left: the histogram of the eigenvalues of $C$. Middle: the loss landscape as a function of $q$, and the global minimum is around the ground truth $q$.  Right: the exact (red) and reconstructed (blue) $\muA$.}
    \label{fig:m1}
  \end{figure}
  
  Figure \ref{fig:m1} summarizes the results. The first plot gives the histogram of the eigenvalues of $C$.  The second plot shows the loss landscape as a function of $q$, and the global minimum is around the ground truth $q$.  The last plot shows the exact (red) and reconstructed (blue) spectral measure of $A$.

\end{example}

\begin{example} The parameters are 
  \begin{itemize}
  \item $N=1024$.
  \item $q = 0.5$.
  \item $\mu_A = \frac{1}{3} \delta_{0.2} + \frac{1}{3} \delta_{0.6} + \frac{1}{3} \delta_{1}$.
  \end{itemize}
  
  \begin{figure}[h!]
    \centering
    \includegraphics[scale=0.27]{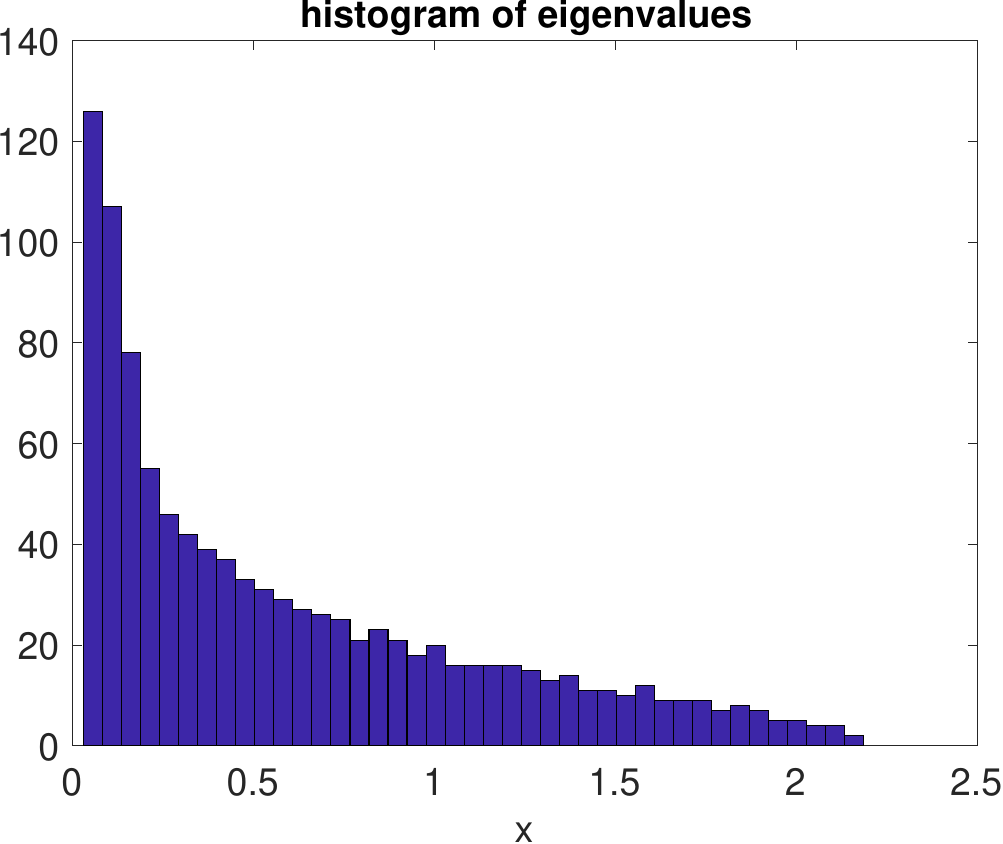}
    \includegraphics[scale=0.27]{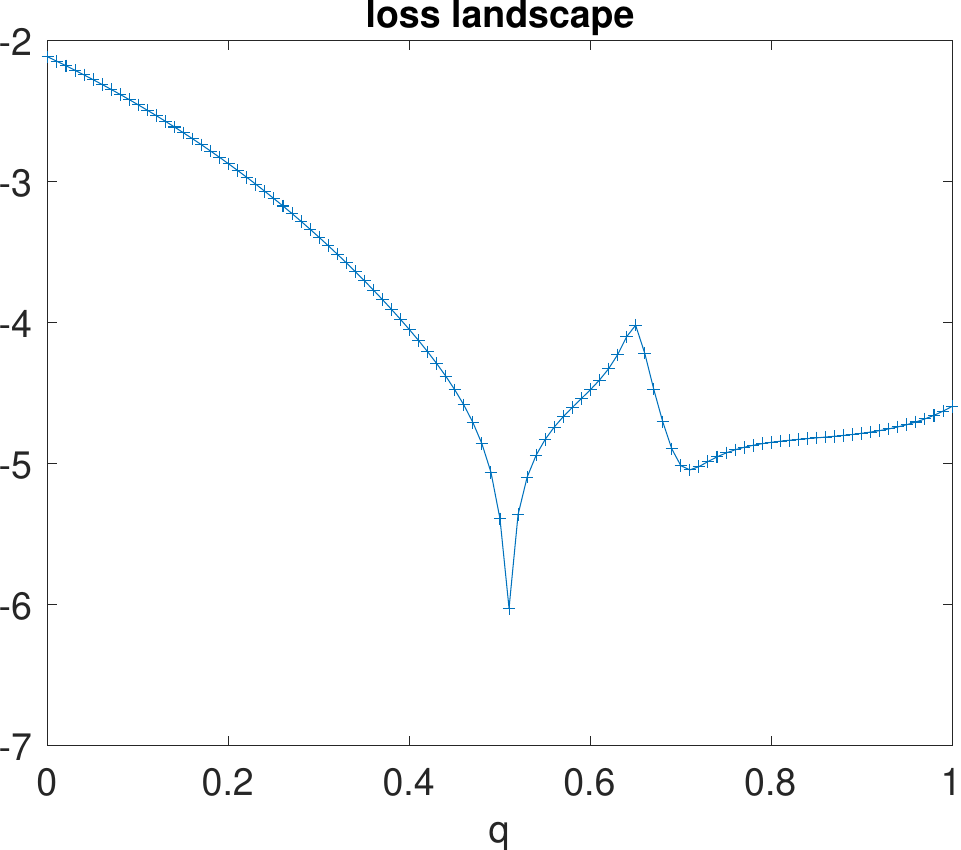}
    \includegraphics[scale=0.27]{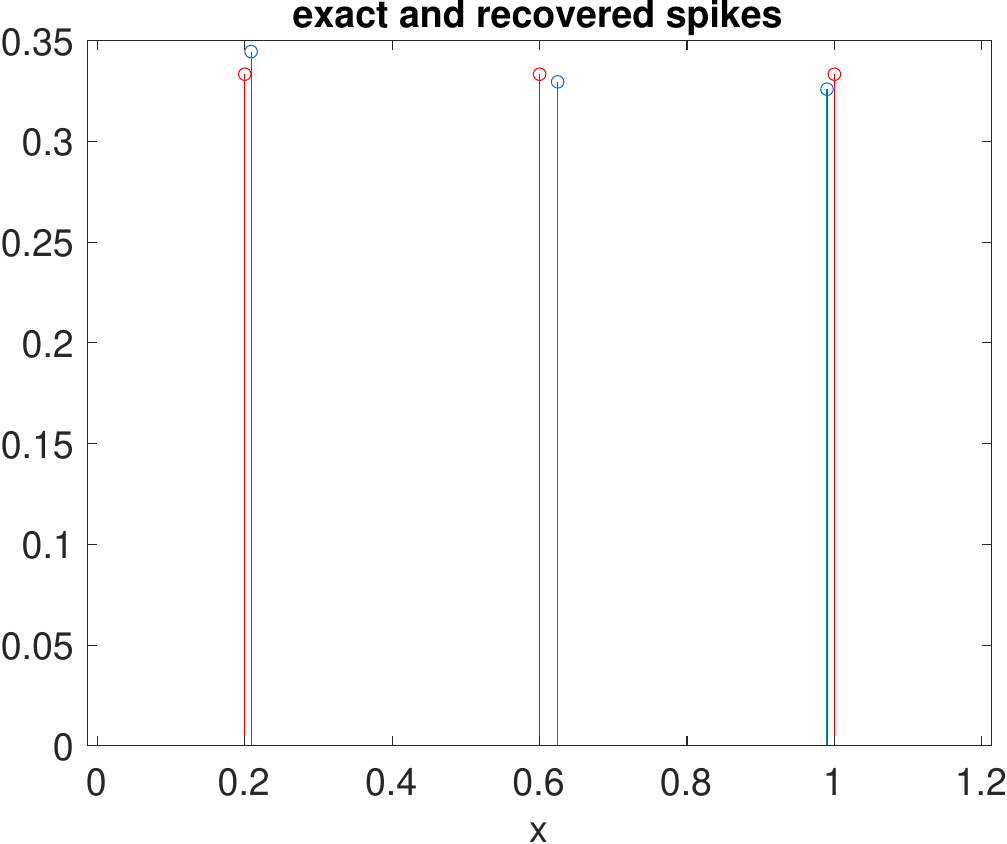}
    \caption{Left: the histogram of the eigenvalues of $C$. Middle: the loss landscape as a function of $q$, and the global minimum is around the ground truth $q$.  Right: the exact (red) and reconstructed (blue) $\muA$.}
    \label{fig:m2}
  \end{figure}

  Figure \ref{fig:m2} summarizes the results. The meanings of the plots are the same as in the previous example.


\end{example}

\begin{example} The parameters are 
  \begin{itemize}
  \item $N=1024$.
  \item $q = 0.75$.
  \item $\mu_A = \frac{1}{3} \delta_{0.2} + \frac{1}{3} \delta_{0.6} + \frac{1}{3} \delta_{1}$.
  \end{itemize}
  
  \begin{figure}[h!]
    \centering
    \includegraphics[scale=0.27]{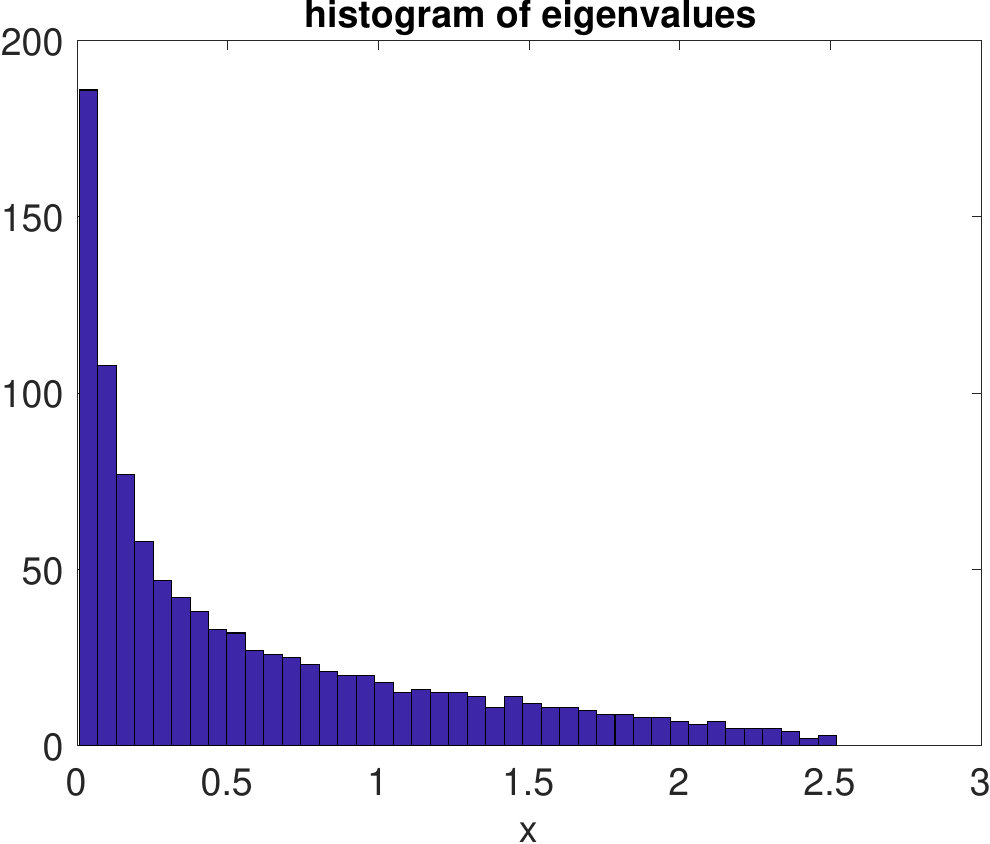}
    \includegraphics[scale=0.27]{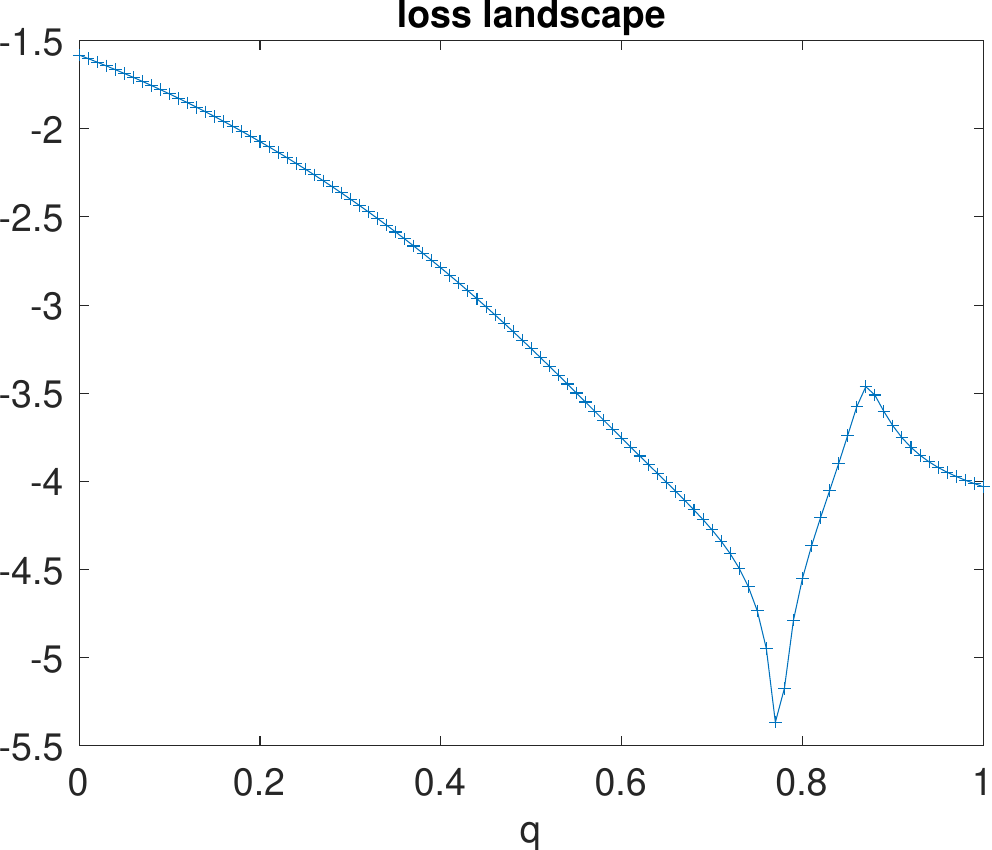}
    \includegraphics[scale=0.27]{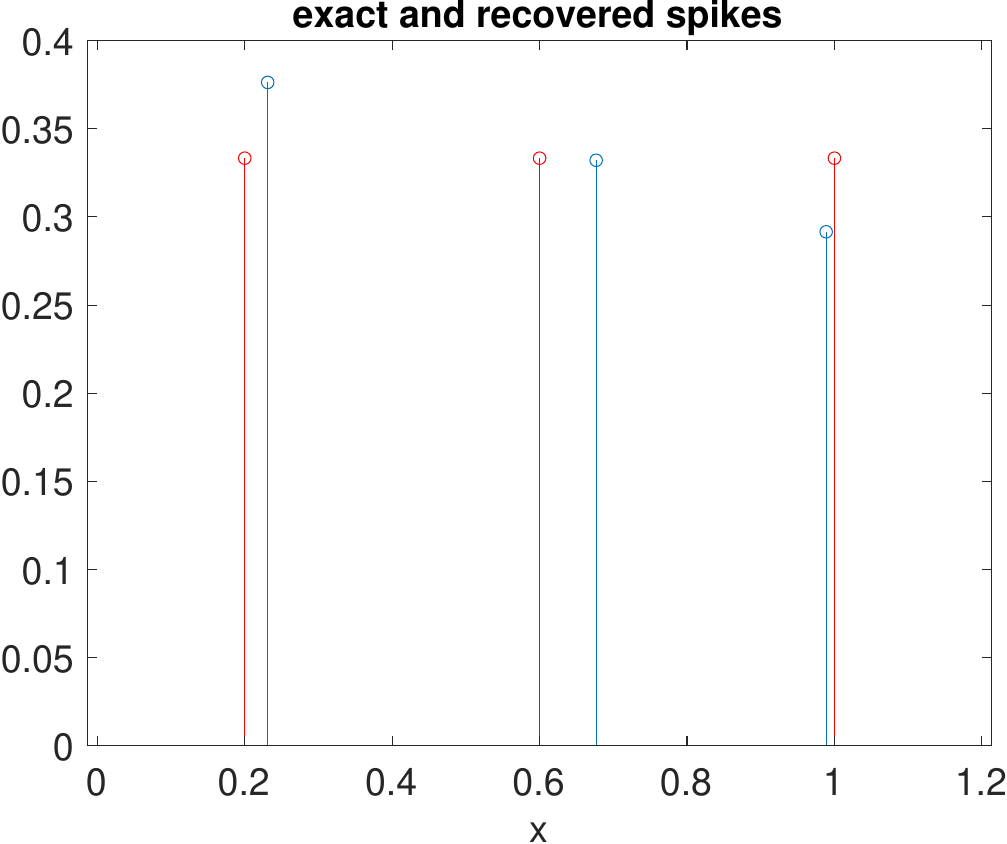}
    \caption{Left: the histogram of the eigenvalues of $C$. Middle: the loss landscape as a function of $q$, and the global minimum is around the ground truth $q$.  Right: the exact (red) and reconstructed (blue) $\muA$.}
    \label{fig:m3}
  \end{figure}

  Figure \ref{fig:m3} summarizes the results. The meanings of the plots are the same as in the previous example.


\end{example}

\section{Discussion}\label{sec:disc}

This section discusses several directions for future work. Here, we consider the most common scenarios for $B$: the Wigner matrix in the additive and the Wishart matrix multiplicative settings. If relevant to applications, other parametric families for $B$ can also be considered.

Second, in the numerical examples, we assume the sparsity $\nx$ is known. As we have seen, discovering the noise level and $\nx$ together robustly from the eigenvalue decay of the $T$ matrix is a non-trivial task. This involves a better understanding of the error of the eigenmatrix method and the finite effect of $N$.

Third, the current approach leverages the asymptotic relationship of the R-transform and the S-transform. For finite $N$, these relationships are approximate for spectral measures. As a result, the current approach has a systematic bias. An important direction is how to incorporate the $N$-dependent corrections for better estimations.

Fourth, the sparsity assumption of the spectral measure of $A$ may not be appropriate for certain applications. For a specific application, if the spectral measure of $A$ arises from some other low-complexity models controlled by a small number of parameters, this approach based on an eigenmatrix might also be useful.

\bibliographystyle{abbrv}

\bibliography{ref}

\end{document}